\documentclass[12pt]{article}%
\usepackage{amsmath}
\usepackage{amsfonts}
\usepackage{amssymb}
\usepackage{graphicx}%
\newtheorem{theorem}{Theorem}
\newtheorem{acknowledgement}[theorem]{Acknowledgement}

\newtheorem{corollary}[theorem]{Corollary}

\newtheorem{lemma}[theorem]{Lemma}

\newtheorem{proposition}[theorem]{Proposition}
\newtheorem{remark}[theorem]{Remark}

\newenvironment{proof}[1][Proof]{\noindent\textbf{#1.} }{\ \rule{0.5em}{0.5em}}
\begin{document}

\title{Asymptotic analysis of generalized Hermite polynomials}
\author{Diego Dominici \thanks{e-mail: dominicd@newpaltz.edu}\\Department of Mathematics\\State University of New York at New Paltz\\75 S. Manheim Blvd. Suite 9\\New Paltz, NY 12561-2443\\USA\\Phone: (845) 257-2607\\Fax: (845) 257-3571 }
\maketitle

\begin{abstract}
We analyze the polynomials $H_{n}^{r}(x)$ considered by Gould and Hopper, which generalize the classical Hermite polynomials. We present the main properties of $H_{n}^{r}(x)$ and derive asymptotic approximations for large values of $n$ from their differential-difference equation, using a discrete ray method. We give numerical examples showing the accuracy of our formulas.

\end{abstract}

Keywords: Hermite polynomials, asymptotic analysis, ray method, differential-difference equations, discrete WKB method.

MSC-class: 33C45 (Primary) 34E05, 34E20 (Secondary)

\section{\strut Introduction}

In \cite{MR0132853} Gould and Hopper generalized the classical Hermite
polynomials by introducing the polynomials $H_{n}^{r}(x,a,p)$ defined by%
\[
H_{n}^{r}(x,a,p)=\left(  -1\right)  ^{n}x^{-a}\exp\left(  px^{r}\right)
\frac{d^{n}}{dx^{n}}\left[  x^{a}\exp\left(  -px^{r}\right)  \right]  .
\]
They derived several properties of $H_{n}^{r}(x,a,p),$ including a generating
function, differentiation, addition and operational formulas. As they
remarked, the case $a=0$ was considered before them by Bell in
\cite{MR1503161}, where he analyzed the polynomials%
\[
\xi_{n}\left(  x,t;r\right)  =\exp\left(  -xt^{r}\right)  \frac{d^{n}}{dt^{n}%
}\exp\left(  xt^{r}\right)  .
\]

Dhawan \cite{MR0252713}, obtained a new generating function, summation
formulas, a hypergeometric representation and some integrals of $H_{n}%
^{r}(x,a,p),$ including%
\[%
{\displaystyle\int\limits_{0}^{\infty}}
x^{a}\exp\left(  -px^{r}\right)  H_{n}^{r}(x,a,p)H_{m}^{s}(x,a,p)dx=0
\]
for $n>ms.$

Kalinowski and Sewery\'{n}ski \cite{MR641527}, \cite{MR1174235}, constructed a
differential equation of order $r$ for $H_{n}^{r}(x,0,p).$ They also proved
that the polynomials $H_{n}^{r}(x,0,p)$ are not orthogonal with respect to the
weight function $\exp\left(  -px^{r}\right)  ,$ $r$ even, in the interval
$\left(  -\infty,\infty\right)  .$

Todorov \cite{MR1206389}, considered the polynomials $H_{n}^{r}(x,0,p)$ in
connection with his analysis of the $n^{th}$ derivative of the composite
function $f\left(  x^{r}\right)  .$ Unaware of Kalinowski and Sewery\'{n}ski's
previous work, he re-derived the differential equation for $H_{n}^{r}(x,0,p).$

Additional properties and further generalizations where studied by Chatterjea
\cite{MR0201691}, Chongdar \cite{MR996881}, Joshi and Prajapat \cite{MR591987}%
, Munot, and Mathur \cite{MR672080}, Rajagopal \cite{MR0104844}, Riordan
\cite{MR0096594}, Saha \cite{MR906780}, Shrivastava \cite{MR0348164} and Singh
and Tiwari \cite{MR1128895}.

In this work, we will consider the polynomials
\begin{equation}
H_{n}^{r}(x)=\left(  -1\right)  ^{n}\exp\left(  x^{r}\right)  \frac{d^{n}%
}{dx^{n}}\exp\left(  -x^{r}\right)  , \label{definition}%
\end{equation}
with $r=2,3,\ldots,\quad n=0,1,\ldots,$ which correspond to the particular
case $H_{n}^{r}(x,0,1).$ \strut Clearly $H_{n}^{2}(x)=$ $H_{n}(x)=$ Hermite
polynomial of degree $n.$ The first few $H_{n}^{r}(x)$ are%
\begin{align*}
H_{0}^{r}(x)  &  =1,\quad H_{1}^{r}(x)=rx^{r-1},\quad H_{2}^{r}(x)=r^{2}%
x^{2\left(  r-1\right)  }-r\left(  r-1\right)  x^{r-2}\\
H_{3}^{r}(x)  &  =r^{3}x^{3\left(  r-1\right)  }-3r^{2}\left(  r-1\right)
x^{2r-3}+r\left(  r-1\right)  \left(  r-2\right)  x^{r-3}%
\end{align*}
and in general%
\[
H_{n}^{r}(x)=r^{n}x^{n\left(  r-1\right)  }-\cdots-\left(  -r\right)
_{n}\ x^{r-n},
\]
where $\left(  \cdot\right)  _{n}$ denotes the Pochhammer symbol. \ Our work
was motivated by the talk "Asymptotics for Hermite type polynomials",
delivered by Professor Wolfgang Gawronski at the conference organized in honor
of the 65th birthday of Nico M. Temme in Santander, Spain on July 4-6, 2005.
In his talk, he considered the behavior of $H_{n}^{r}(x)$ and its zeros using
Plancherel-Rotach type asymptotics.

\strut In Section 2, we present the basic properties of $H_{n}^{r}(x).$
Although some are not new, we present proofs of all of them for completion
purposes. In Section 3, we present the asymptotic analysis using a modified
ray method, developed by Dosdale, Duggan and Morgan in \cite{MR0361328} and
formalized by Costin and Costin in \cite{MR1373150}. In a previous work
\cite{hermitedif}, we successfully applied the same technique to the classical
Hermite polynomials $H_{n}^{2}(x).$ Section 4 contains our main result and
supporting numerical examples.

\section{\strut Properties}

The Hermite polynomials admit the simple hypergeometric representation
\cite{koekoek94askeyscheme}%
\[
H_{n}(x)=\left(  2x\right)  ^{n}\ _{2}F_{0}\left(  \left.
\begin{array}
[c]{c}%
-\frac{n}{2},-\frac{\left(  n-1\right)  }{2}\\
\_
\end{array}
\right\vert -x^{-2}\right)  ,
\]
where $_{p}F_{q}\left[  \cdot\right]  $ is the hypergeometric function
\cite{MR2128719}. On the other hand, for $H_{n}^{r}(x)$ we need to consider
the extension of $_{p}F_{q}$ given by Meijer's $G$-function \cite{MR698779}.

\begin{proposition}
The polynomials $H_{n}^{r}(x)$ can be represented in terms of Meijer's
$G$-function by%
\[
H_{n}^{r}(x)=\exp\left(  x^{r}\right)  \left(  -\frac{r}{x}\right)
^{n}G_{r,r+1}^{1,r}\left(  x^{r}\left\vert
\genfrac{}{}{0pt}{}{1-\frac{1}{r},1-\frac{2}{r},\cdots,0}{0,1+\frac{n-1}%
{r},1+\frac{n-2}{r},\cdots,\frac{n}{r}}%
\right.  \right)  .
\]

\end{proposition}

\begin{proof}
We have%
\begin{equation}
\frac{d^{n}}{dx^{n}}\exp\left(  -x^{r}\right)  =\sum\limits_{j=0}^{\infty
}\frac{d^{n}}{dx^{n}}\frac{\left(  -x^{r}\right)  ^{j}}{j!}=\frac{1}{x^{n}%
}\sum\limits_{j=0}^{\infty}\frac{\left(  rj\right)  !}{\left(  rj-n\right)
!}\frac{\left(  -x^{r}\right)  ^{j}}{j!}.\label{diffexp}%
\end{equation}
Using the multiplication formula for the Gamma function \cite{MR1376370}, we
obtain%
\[
\frac{\left(  rj\right)  !}{\left(  rj-n\right)  !}=\frac{\Gamma\left(
rj+1\right)  }{\Gamma\left(  rj-n+1\right)  }=\frac{r^{rj+1}%
{\displaystyle\prod\limits_{i=0}^{r-1}}
\Gamma\left(  \frac{i+1}{r}+j\right)  }{r^{rj-n+1}%
{\displaystyle\prod\limits_{i=0}^{r-1}}
\Gamma\left(  \frac{i+1-n}{r}+j\right)  }%
\]%
\begin{equation}
=r^{n}%
{\displaystyle\prod\limits_{i=1}^{r}}
\frac{\Gamma\left(  \frac{i}{r}+j\right)  }{\Gamma\left(  \frac{i-n}%
{r}+j\right)  }=r^{n}%
{\displaystyle\prod\limits_{i=1}^{r}}
\frac{\Gamma\left(  \frac{i}{r}\right)  \left(  \frac{i}{r}\right)  _{j}%
}{\Gamma\left(  \frac{i-n}{r}\right)  \left(  \frac{i-n}{r}\right)  _{j}%
}.\label{Gamma1}%
\end{equation}
Replacing (\ref{Gamma1}) in (\ref{diffexp}) we get%
\begin{gather*}
\frac{d^{n}}{dx^{n}}\exp\left(  -x^{r}\right)  =\left(  \frac{r}{x}\right)
^{n}%
{\displaystyle\prod\limits_{i=1}^{r}}
\frac{\Gamma\left(  \frac{i}{r}\right)  }{\Gamma\left(  \frac{i-n}{r}\right)
}\sum\limits_{j=0}^{\infty}%
{\displaystyle\prod\limits_{i=1}^{r}}
\frac{\left(  \frac{i}{r}\right)  _{j}}{\left(  \frac{i-n}{r}\right)  _{j}%
}\frac{\left(  -x^{r}\right)  ^{j}}{j!}\\
=\left(  \frac{r}{x}\right)  ^{n}%
{\displaystyle\prod\limits_{i=1}^{r}}
\frac{\Gamma\left(  \frac{i}{r}\right)  }{\Gamma\left(  \frac{i-n}{r}\right)
}\ _{r}F_{_{_{r}}}\left[  \left.
\genfrac{}{}{0pt}{}{\frac{1}{r},\frac{2}{r},\cdots,1}{\frac{1-n}{r},\frac
{2-n}{r},\cdots,\frac{r-n}{r}}%
\right\vert -x^{r}\right]  .
\end{gather*}
Thus, from the definition of $H_{n}^{r}(x),$ it follows that%
\begin{equation}
H_{n}^{r}(x)=\exp\left(  x^{r}\right)  \left(  -\frac{r}{x}\right)  ^{n}%
{\displaystyle\prod\limits_{i=1}^{r}}
\frac{\Gamma\left(  \frac{i}{r}\right)  }{\Gamma\left(  \frac{i-n}{r}\right)
}\ _{r}F_{_{_{r}}}\left[  \left.
\genfrac{}{}{0pt}{}{\frac{1}{r},\frac{2}{r},\cdots,1}{\frac{1-n}{r},\frac
{2-n}{r},\cdots,\frac{r-n}{r}}%
\right\vert -x^{r}\right]  .\label{hyperg}%
\end{equation}

Using the relation between the hypergeometric function and Meijer's
$G$-function \cite{MR1773820}%
\[
\frac{%
{\displaystyle\prod\limits_{i=1}^{p}}
\Gamma\left(  a_{i}\right)  }{%
{\displaystyle\prod\limits_{i=1}^{q}}
\Gamma\left(  b_{i}\right)  }\ _{p}F_{q}\left[  \left.
\genfrac{}{}{0pt}{}{a_{1},\cdots a_{p}}{b_{1},\cdots b_{q}}%
\right\vert x\right]  =G_{p,q+1}^{1,p}\left(  -x\left\vert \left.
\genfrac{}{}{0pt}{}{1-a_{1},\cdots1-a_{p}}{0,1-b_{1},\cdots1-b_{q}}%
\right\vert \right.  \right)
\]
in (\ref{hyperg}), the results follow.
\end{proof}

\begin{proposition}
The polynomials $H_{n}^{r}(x)$ satisfy the differential-difference equation%
\begin{equation}
H_{n+1}^{r}(x)+\frac{d}{dx}H_{n}^{r}(x)=rx^{r-1}H_{n}^{r}(x).
\label{Hdiffdiff}%
\end{equation}

\end{proposition}

\begin{proof}
The result follows immediately from the Rodrigues formula (\ref{definition}),
since%
\begin{gather*}
\frac{d}{dx}H_{n}^{r}(x)=\left(  -1\right)  ^{n}rx^{r-1}\exp\left(
x^{r}\right)  \frac{d^{n}}{dx^{n}}\exp\left(  -x^{r}\right) \\
+\left(  -1\right)  ^{n}\exp\left(  x^{r}\right)  \frac{d^{n+1}}{dx^{n+1}}%
\exp\left(  -x^{r}\right)  =rx^{r-1}H_{n}^{r}(x)-H_{n+1}^{r}(x).
\end{gather*}

\end{proof}

When $r=2,$ we recover the well-known formula for the Hermite polynomials
\cite{MR0350075}%
\[
H_{n+1}(x)+H_{n}^{\prime}(x)=2xH_{n}(x).
\]

\begin{proposition}
The polynomials $H_{n}^{r}(x)$ have the exponential generating function%
\begin{equation}
G(x,t)=%
{\displaystyle\sum\limits_{n=0}^{\infty}}
H_{n}^{r}(x)\frac{t^{n}}{n!}=\exp\left[  x^{r}-(x-t)^{r}\right]  .
\label{Generating}%
\end{equation}

\end{proposition}

\begin{proof}
From (\ref{definition}) we get
\begin{align*}
G(x,t)  &  =\exp\left(  x^{r}\right)
{\displaystyle\sum\limits_{n=0}^{\infty}}
\left(  -1\right)  ^{n}\frac{t^{n}}{n!}\left[  \frac{d^{n}}{du^{n}}e^{-u^{r}%
}\right]  _{u=x}\\
&  =\exp\left(  x^{r}\right)
{\displaystyle\sum\limits_{n=0}^{\infty}}
\frac{t^{n}}{n!}\left[  \frac{d^{n}}{dt^{n}}e^{-\left(  x-t\right)  ^{r}%
}\right]  _{t=0}\\
&  =\exp\left(  x^{r}\right)  \exp\left[  -\left(  x-t\right)  ^{r}\right]  .
\end{align*}

\end{proof}

In particular, for $r=2,$ we have%
\[
\exp\left[  x^{2}-(x-t)^{2}\right]  =\exp\left(  2tx-t^{2}\right)  ,
\]
which is the exponential generating function of the Hermite polynomials
\cite{MR0350075}.

\begin{proposition}
The polynomials $H_{n}^{r}(x)$ admit the explicit representation%
\begin{equation}
H_{n}^{r}(x)=%
{\displaystyle\sum\limits_{k=\left\lfloor \frac{n}{r}\right\rfloor }^{n}}
C_{k}^{n}(r)x^{rk-n}, \label{Hrepresent}%
\end{equation}
where%
\begin{equation}
C_{k}^{n}(r)=\frac{\left(  -1\right)  ^{n}n!}{k!}%
{\displaystyle\sum\limits_{j=0}^{k}}
\left(  -1\right)  ^{j}\dbinom{k}{j}\dbinom{rj}{n}. \label{Cn}%
\end{equation}

\end{proposition}

\begin{proof}
From (\ref{definition}) we have%
\begin{align*}
H_{n}^{r}(x) &  =\left(  -1\right)  ^{n}\exp\left(  x^{r}\right)  \frac{d^{n}%
}{dx^{n}}\exp\left(  -x^{r}\right)  \\
&  =\left(  -1\right)  ^{n}\left[  \sum\limits_{k=0}^{\infty}\frac{\left(
x^{r}\right)  ^{k}}{k!}\right]  \left[  \sum\limits_{j=0}^{\infty}\frac{d^{n}%
}{dx^{n}}\frac{\left(  -x^{r}\right)  ^{j}}{j!}\right]  \\
&  =\frac{\left(  -1\right)  ^{n}n!}{x^{n}}\left[  \sum\limits_{k=0}^{\infty
}\frac{\left(  x^{r}\right)  ^{k}}{k!}\right]  \left[  \sum\limits_{j=0}%
^{\infty}\frac{\left(  -1\right)  ^{j}}{j!}\dbinom{rj}{n}\left(  x^{r}\right)
^{j}\right]  ,
\end{align*}
and hence \cite[(2.8)]{MR2172781}%
\[
H_{n}^{r}(x)=\frac{\left(  -1\right)  ^{n}n!}{x^{n}}\sum\limits_{k=0}^{\infty
}\left[  \sum\limits_{j=0}^{k}\frac{\left(  -1\right)  ^{j}}{j!}\dbinom{rj}%
{n}\frac{1}{\left(  k-j\right)  !}\right]  \left(  x^{r}\right)  ^{k}.
\]
If $rk<n,$ we have
\[
\dbinom{rj}{n}=0,\quad0\leq j\leq k
\]
and therefore,%
\[
H_{n}^{r}(x)=\frac{\left(  -1\right)  ^{n}n!}{x^{n}}\sum
\limits_{k=\left\lfloor \frac{n}{r}\right\rfloor }^{\infty}\left[
\sum\limits_{j=0}^{k}\frac{\left(  -1\right)  ^{j}}{j!}\dbinom{rj}{n}\frac
{1}{\left(  k-j\right)  !}\right]  \left(  x^{r}\right)  ^{k},
\]
from which (\ref{Hrepresent}) follows.
\end{proof}

Note that when $r=2,$ (\ref{Hrepresent}) reduces to the well-known
representation of the Hermite polynomials \cite{MR0350075}%
\begin{gather*}
H_{n}(x)=%
{\displaystyle\sum\limits_{k=\left\lfloor \frac{n}{2}\right\rfloor }^{n}}
\frac{\left(  -1\right)  ^{n}n!}{k!}\left[
{\displaystyle\sum\limits_{j=0}^{k}}
\left(  -1\right)  ^{j}\dbinom{k}{j}\dbinom{2j}{n}\right]  x^{2k-n}\\
=%
{\displaystyle\sum\limits_{k=\left\lfloor \frac{n}{2}\right\rfloor }^{n}}
\left(  -1\right)  ^{n}n!\frac{\left(  -1\right)  ^{k}2^{2k-n}}{\left(
2k-n\right)  !(n-k)!}x^{2k-n}=%
{\displaystyle\sum\limits_{k=0}^{\left\lfloor \frac{n}{2}\right\rfloor }}
\frac{\left(  -1\right)  ^{k}n!}{\left(  n-2k\right)  !k!}\left(  2x\right)
^{n-2k}.
\end{gather*}

\begin{corollary}
\label{symm}The polynomials $H_{n}^{r}(x)$ satisfy%
\begin{equation}
H_{n}^{r}(x)=\omega^{n}H_{n}^{r}(\omega x),\quad\text{with \ }\omega^{r}=1.
\label{symmetry}%
\end{equation}

\end{corollary}

The symmetry relation (\ref{symmetry}) generalizes the reflection formula of
the Hermite polynomials \cite{MR0350075}%
\[
H_{n}(x)=\left(  -1\right)  ^{n}H_{n}(-x).
\]

\begin{figure}[ptb]
\begin{center}
\rotatebox{270} {\resizebox{!}{5in}{\includegraphics{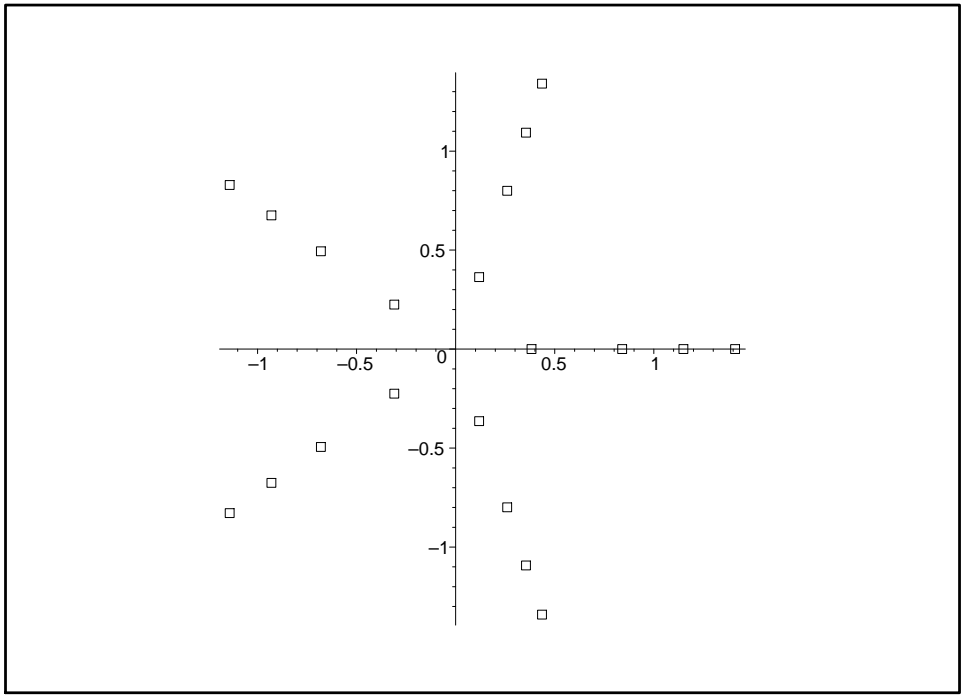}}}
\end{center}
\caption{A plot of the zeros of $H_{5}^{5}(x)$ in the complex plane.}%
\label{zeros}%
\end{figure}

\begin{remark}
It follows from (\ref{symmetry}) that it is enough to analyze $H_{n}^{r}(x)$
in the region $\left\vert \arg(x)\right\vert \leq$ $\frac{\pi}{r}$ of the
complex plane. In particular, the roots of $H_{n}^{r}(x)$ are completely
determined by its positive real roots. To show this, in Figure \ref{zeros} we
plot the roots of $H_{5}^{5}(x)$ in the complex plane.
\end{remark}

\begin{proposition}
The polynomials $H_{n}^{r}(x)$ can be represented by%
\[
H_{n}^{r}(x)=\left(  -1\right)  ^{n}n!%
{\displaystyle\sum\limits_{\substack{N_{1}+\cdots+N_{r}=n\\n\geq0}}}
{\displaystyle\prod\limits_{j=1}^{r}}
\Omega_{j}\left(  N_{j}\right)  \binom{r}{j}^{\frac{N_{j}}{j}}\frac{\left(
-1\right)  ^{\frac{N_{j}}{j}}}{\left(  \frac{N_{j}}{j}\right)  !}%
x^{\frac{\left(  r-j\right)  }{j}N_{j}}.
\]
In particular, we have%
\begin{equation}
H_{n}^{r}(0)=\left(  -1\right)  ^{\frac{r-1}{r}n}\frac{n!}{\left(  \frac{n}%
{r}\right)  !}\Omega_{r}\left(  n\right)  , \label{H0}%
\end{equation}
where \cite[(2.32)]{MR2172781}%
\begin{equation}
\Omega_{r}\left(  k\right)  =\frac{1}{r}%
{\displaystyle\sum\limits_{j=1}^{r}}
\exp\left(  \frac{2\pi jk\mathrm{i}}{r}\right)  =\left\{
\begin{array}
[c]{c}%
1\quad\quad\text{if}\quad\quad r|_{k}\\
0\quad\text{otherwise}%
\end{array}
\right.  . \label{Omega}%
\end{equation}

\end{proposition}

\begin{proof}
From (\ref{Generating}) we get%
\begin{equation}%
{\displaystyle\sum\limits_{n=0}^{\infty}}
\frac{H_{n}^{r}(x)}{n!}t^{n} \label{part1}%
\end{equation}%
\begin{gather*}
=\exp\left[  x^{r}-(x-t)^{r}\right]  =\exp\left[  x^{r}-%
{\displaystyle\sum\limits_{j=0}^{r}}
\binom{r}{j}\left(  -1\right)  ^{j}x^{r-j}t^{j}\right] \\
=\exp\left[
{\displaystyle\sum\limits_{j=1}^{r}}
\binom{r}{j}\left(  -1\right)  ^{j+1}x^{r-j}t^{j}\right]  =%
{\displaystyle\prod\limits_{j=1}^{r}}
\exp\left[  \binom{r}{j}\left(  -1\right)  ^{j+1}x^{r-j}t^{j}\right]
\end{gather*}%
\begin{align*}
&  =%
{\displaystyle\prod\limits_{j=1}^{r}}
{\displaystyle\sum\limits_{N_{j}=0}^{\infty}}
\binom{r}{j}^{N_{j}}\frac{\left(  -1\right)  ^{\left(  j+1\right)  N_{j}}%
}{\left(  N_{j}\right)  !}x^{\left(  r-j\right)  N_{j}}t^{jN_{j}}\\
&  =%
{\displaystyle\prod\limits_{j=1}^{r}}
{\displaystyle\sum\limits_{N_{j}=0}^{\infty}}
\Omega_{j}\left(  N_{j}\right)  \binom{r}{j}^{\frac{N_{j}}{j}}\frac{\left(
-1\right)  ^{\frac{\left(  j-1\right)  }{j}N_{j}}}{\left(  \frac{N_{j}}%
{j}\right)  !}x^{\frac{\left(  r-j\right)  }{j}N_{j}}t^{N_{j}}%
\end{align*}%
\begin{equation}
=%
{\displaystyle\sum\limits_{\substack{N_{1}+\cdots+N_{r}=n\\n\geq0}}}
\left[
{\displaystyle\prod\limits_{j=1}^{r}}
\Omega_{j}\left(  N_{j}\right)  \binom{r}{j}^{\frac{N_{j}}{j}}\frac{\left(
-1\right)  ^{\frac{\left(  j-1\right)  }{j}N_{j}}}{\left(  \frac{N_{j}}%
{j}\right)  !}x^{\frac{\left(  r-j\right)  }{j}N_{j}}\right]  t^{n},
\label{part2}%
\end{equation}
where we have used equation (2.10) in \cite{MR2172781}.

Comparing (\ref{part1}) and (\ref{part2}), we conclude that%
\begin{align*}
H_{n}^{r}(x)  &  =n!%
{\displaystyle\sum\limits_{\substack{N_{1}+\cdots+N_{r}=n\\n\geq0}}}
{\displaystyle\prod\limits_{j=1}^{r}}
\Omega_{j}\left(  N_{j}\right)  \binom{r}{j}^{\frac{N_{j}}{j}}\frac{\left(
-1\right)  ^{\frac{\left(  j+1\right)  }{j}N_{j}}}{\left(  \frac{N_{j}}%
{j}\right)  !}x^{\frac{\left(  r-j\right)  }{j}N_{j}}\\
&  =\left(  -1\right)  ^{n}n!%
{\displaystyle\sum\limits_{\substack{N_{1}+\cdots+N_{r}=n\\n\geq0}}}
{\displaystyle\prod\limits_{j=1}^{r}}
\Omega_{j}\left(  N_{j}\right)  \binom{r}{j}^{\frac{N_{j}}{j}}\frac{\left(
-1\right)  ^{\frac{N_{j}}{j}}}{\left(  \frac{N_{j}}{j}\right)  !}%
x^{\frac{\left(  r-j\right)  }{j}N_{j}}%
\end{align*}
In particular, when $x=0,$ we must have
\[
N_{1}=\cdots=N_{r-1}=0\text{ \ \ and \ \ }N_{r}=n,
\]
from which (\ref{H0}) follows.
\end{proof}

\begin{proposition}
The polynomials $H_{n}^{r}(x)$ satisfy the recurrence relation%
\begin{equation}
H_{n+1}^{r}(x)=r%
{\displaystyle\sum_{k=0}^{r-1}}
\left(  -1\right)  ^{k}k!\binom{n}{k}\dbinom{r-1}{k}x^{r-1-k}H_{n-k}^{r}(x).
\label{Hrecurr}%
\end{equation}

\end{proposition}

\begin{proof}
From (\ref{Generating}), we have%
\[
\frac{\partial G}{\partial t}=r\left(  x-t\right)  ^{r-1}G
\]
or [Wilf],%
\[%
{\displaystyle\sum\limits_{n=0}^{\infty}}
H_{n+1}^{r}(x)\frac{t^{n}}{n!}=r\left[
{\displaystyle\sum\limits_{k=0}^{\infty}}
\left(  -1\right)  ^{k}\binom{r-1}{k}x^{r-1-k}t^{k}\right]  \left[
{\displaystyle\sum\limits_{n=0}^{\infty}}
H_{n}^{r}(x)\frac{t^{n}}{n!}\right]  .
\]
Comparing coefficients of $t$ we get%
\[
\frac{H_{n+1}^{r}(x)}{n!}=r%
{\displaystyle\sum\limits_{k=0}^{n}}
\left(  -1\right)  ^{k}\binom{r-1}{k}x^{r-1-k}\frac{H_{n-k}^{r}(x)}{\left(
n-k\right)  !}%
\]
and (\ref{Hrecurr}) follows.
\end{proof}

When $r=2,$ (\ref{Hrecurr}) reduces to the three-term recurrence relation for
the Hermite polynomials \cite{MR0350075}%
\[
H_{n+1}(x)=2xH_{n}(x)-2nH_{n-1}(x).
\]

\begin{proposition}
Let $u_{n}(x)=\exp\left(  -x^{r}\right)  H_{n}^{r}(x).$ Then,%
\begin{equation}
u_{n}^{\left(  r\right)  }+r%
{\textstyle\sum\limits_{k=0}^{r-1}}
\dbinom{r-1}{k}\frac{\left(  n+r-1\right)  !}{\left(  n+k\right)  !}x^{k}%
u_{n}^{\left(  k\right)  }=0, \label{Hdiffeq}%
\end{equation}
where
\[
u_{n}^{\left(  k\right)  }=\frac{d^{k}}{dx^{k}}u_{n}\left(  x\right)  .
\]

\end{proposition}

\begin{proof}
It is clear from (\ref{definition}) that%
\begin{equation}
u_{n+k}(x)=\left(  -1\right)  ^{k}u_{n}^{\left(  k\right)  }. \label{u}%
\end{equation}
From (\ref{Hrecurr}), we get%
\begin{align*}
u_{n+r}(x)  &  =r%
{\displaystyle\sum_{k=0}^{r-1}}
\left(  -1\right)  ^{k}k!\binom{n+r-1}{k}\dbinom{r-1}{k}x^{r-1-k}%
u_{n+r-1-k}(x)\\
&  =r%
{\displaystyle\sum_{k=0}^{r-1}}
\left(  -1\right)  ^{r-1-k}\left(  r-1-k\right)  !\binom{n+r-1}{r-1-k}%
\dbinom{r-1}{r-1-k}x^{k}u_{n+k}(x).
\end{align*}
Thus,%
\begin{equation}
u_{n+r}(x)=r%
{\displaystyle\sum_{k=0}^{r-1}}
\left(  -1\right)  ^{r-1-k}\dbinom{r-1}{k}\frac{\left(  n+r-1\right)
!}{\left(  n+k\right)  !}x^{k}u_{n+k}(x). \label{diff1}%
\end{equation}
Using (\ref{u}) in (\ref{diff1}), the result follows.
\end{proof}

For the case $r=2,$ (\ref{Hdiffeq}) gives%
\[
\left[  \exp\left(  -x^{2}\right)  H_{n}\right]  ^{\prime\prime}+2x\left[
\exp\left(  -x^{2}\right)  H_{n}\right]  ^{\prime}+2\left(  n+1\right)
\exp\left(  -x^{2}\right)  H_{n}=0,
\]
which is equivalent to the differential equation of the Hermite polynomials
\cite{MR0350075}%
\[
H_{n}^{\prime\prime}-2xH_{n}^{\prime}+2nH_{n}=0.
\]

\section{Asymptotic analysis}

We seek an approximative solution for (\ref{Hdiffdiff}) of the form
\begin{equation}
H_{n}^{r}(x)\sim\exp\left[  f(x,n)+g(x,n)\right]  ,\quad n\rightarrow
\infty\label{anszat}%
\end{equation}
with%
\[
g=o(f),\quad n\rightarrow\infty.
\]
Since $H_{0}^{r}(x)=1,$ we must have%
\begin{equation}
f(x,0)=0\text{ \ \ and \ \ }g(x,0)=0. \label{initial}%
\end{equation}
Using (\ref{anszat}) in (\ref{Hdiffdiff}), we have%
\begin{gather}
\exp\left(  f+\frac{\partial f}{\partial n}+\frac{1}{2}\frac{\partial^{2}%
f}{\partial n^{2}}+g+\frac{\partial g}{\partial n}\right) \label{asymp1}\\
+\left(  \frac{\partial f}{\partial x}+\frac{\partial g}{\partial x}\right)
\exp\left(  f+g\right)  =rx^{r-1}\exp\left(  f+g\right)  ,\nonumber
\end{gather}
where we have used%
\[
f(x,n+1)=f(x,n)+\frac{\partial f}{\partial n}(x,n)+\frac{1}{2}\frac
{\partial^{2}f}{\partial n^{2}}(x,n)+\cdots.
\]

From (\ref{asymp1}) we obtain, to leading order, the \textit{eikonal} equation%
\begin{equation}
\exp\left(  \frac{\partial f}{\partial n}\right)  +\frac{\partial f}{\partial
x}-rx^{r-1}=0, \label{eikonal}%
\end{equation}
and
\[
\exp\left(  \frac{1}{2}\frac{\partial^{2}f}{\partial n^{2}}+\frac{\partial
g}{\partial n}\right)  +\frac{\partial g}{\partial x}\exp\left(
-\frac{\partial f}{\partial n}\right)  -1=0,
\]
or, to leading order, the \textit{transport }equation%
\begin{equation}
\frac{1}{2}\frac{\partial^{2}f}{\partial n^{2}}+\frac{\partial g}{\partial
n}+\frac{\partial g}{\partial x}\exp\left(  -\frac{\partial f}{\partial
n}\right)  =0. \label{transport}%
\end{equation}

\subsection{The rays}

To solve (\ref{eikonal}), we use the method of characteristics, which we
briefly review. Given the first order partial differential equation%
\[
F\left(  x,n,f,p,q\right)  =0,\text{ \ \ with \ \ }\ p=\frac{\partial
f}{\partial x},\quad q=\frac{\partial f}{\partial n},
\]
we search for a solution \ $f(x,n)$ by solving the system of \textquotedblleft
characteristic equations\textquotedblright\
\begin{align*}
\frac{dx}{dt}  &  =\frac{\partial F}{\partial p},\quad\frac{dn}{dt}%
=\frac{\partial F}{\partial q},\\
\frac{dp}{dt}  &  =-\frac{\partial F}{\partial x}-p\frac{\partial F}{\partial
f},\quad\frac{dq}{dt}=-\frac{\partial F}{\partial n}-q\frac{\partial
F}{\partial f},\\
\frac{df}{dt}  &  =p\frac{\partial F}{\partial p}+q\frac{\partial F}{\partial
q},
\end{align*}
with initial conditions%
\begin{equation}
F\left[  x(0,s),n(0,s),f(0,s),p(0,s),q(0,s)\right]  =0, \label{initial1}%
\end{equation}
and%
\begin{equation}
\quad\frac{d}{ds}f(0,s)=p(0,s)\frac{d}{ds}x(0,s)+q(0,s)\frac{d}{ds}n(0,s),
\label{initial2}%
\end{equation}
where we now consider $\left\{  x,n,f,p,q\right\}  $ to all be functions of
the variables $t$ and $s.$

For the eikonal equation (\ref{eikonal}), we have%
\begin{equation}
F\left(  x,n,f,p,q\right)  =e^{q}+p-rx^{r-1} \label{eikonal1}%
\end{equation}
and therefore the characteristic equations are%
\begin{equation}
\frac{dx}{dt}=1,\quad\frac{dn}{dt}=e^{q},\quad\frac{dp}{dt}=r\left(
r-1\right)  x^{r-2},\quad\frac{dq}{dt}=0, \label{charac1}%
\end{equation}
and%
\begin{equation}
\frac{df}{dt}=p+qe^{q}. \label{eqf}%
\end{equation}
Solving (\ref{charac1}) subject to the initial conditions%
\[
x(0,s)=s,\quad n(0,s)=0,\quad q\left(  0,s\right)  =A(s),
\]
with $A(s)$ to be determined, we obtain%
\[
x=t+s,\quad n=te^{A},\quad p=r\left(  t+s\right)  ^{r-1}+B(s),\quad q=A,
\]
for some function $B(s).$ From (\ref{initial1}) we have%
\[
e^{A}+rs^{r-1}+B-rs^{r-1}=0
\]
and $B=-e^{A}.$ Thus,%
\begin{equation}
x=t+s,\quad n=te^{A},\quad p=r\left(  t+s\right)  ^{r-1}-e^{A},\quad q=A.
\label{char1}%
\end{equation}

Since (\ref{initial}) implies that $f(0,s)=0,$ we have from (\ref{initial2})
and (\ref{char1})%
\[
\left(  rs^{r-1}-e^{A}\right)  \times1+A\times0=0.
\]
Hence, $A=\ln\left(  rs^{r-1}\right)  $ and therefore%
\begin{equation}
x=t+s,\quad n=rs^{r-1}t, \label{rays}%
\end{equation}%
\begin{equation}
\quad p=r\left[  \left(  t+s\right)  ^{r-1}-s^{r-1}\right]  ,\quad
q=\ln\left(  rs^{r-1}\right)  , \label{pq}%
\end{equation}
with $s>0.$

\subsection{The caustic}

\begin{figure}[ptb]
\begin{center}
\rotatebox{270} {\resizebox{!}{5in}{\includegraphics{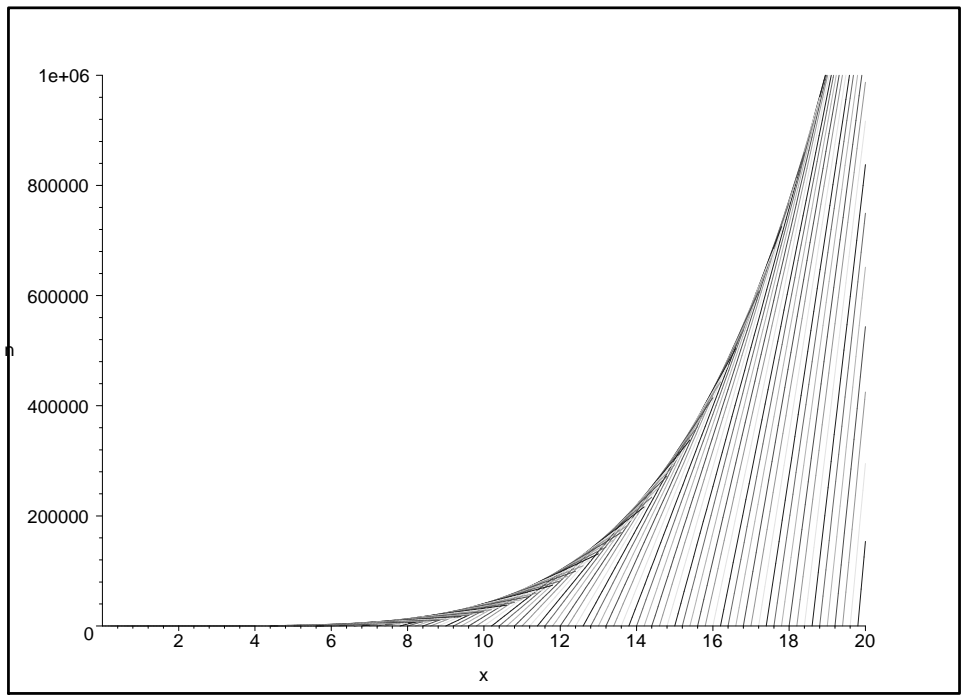}}}
\end{center}
\caption{A sketch of the rays with $r=5$ showing the caustic.}%
\label{ray}%
\end{figure}

Sketching the rays (\ref{rays}), we observe that they fill the region
$x>X_{c}\left(  n\right)  $, where $X_{c}\left(  n\right)  $ is the
\textit{caustic,} i.e., the points in the $(x,n)$-plane at which the Jacobian%
\begin{equation}
J\left(  t,s\right)  =\frac{\partial x}{\partial t}\frac{\partial n}{\partial
s}-\frac{\partial x}{\partial s}\frac{\partial n}{\partial t}=rs^{r-2}\left[
(r-1)t-s\right]  \label{J}%
\end{equation}
is zero (see Fig \ref{ray}).

From (\ref{J}) we have, for $s>0$,
\begin{equation}
J(t,s)=0\Leftrightarrow s=(r-1)t. \label{Js}%
\end{equation}
Using (\ref{Js}) in (\ref{rays}) we obtain
\begin{equation}
J(t,s)=0\Leftrightarrow t=\frac{x}{r}\Leftrightarrow s=\lambda x \label{J1}%
\end{equation}
with%
\begin{equation}
\lambda=\frac{r-1}{r},\quad\frac{1}{2}\leq\lambda<1. \label{lambda}%
\end{equation}
From (\ref{Js}) and (\ref{J1}) we conclude that
\begin{equation}
X_{c}\left(  n\right)  =\lambda^{-\lambda}n^{1-\lambda}. \label{Xc}%
\end{equation}

\subsection{$\allowbreak$The functions $f$ and $g$}

Using (\ref{pq}) in (\ref{eqf}) we have%
\begin{equation}
\frac{df}{dt}=r\left[  \left(  t+s\right)  ^{r-1}-s^{r-1}\right]  +\ln\left(
rs^{r-1}\right)  rs^{r-1}, \label{eqf1}%
\end{equation}
while (\ref{initial}) implies that $f(0,s)=0.$ Solving (\ref{eqf1}), we obtain%
\begin{equation}
f(t,s)=\left(  t+s\right)  ^{r}-s^{r}+\left[  \ln\left(  rs^{r-1}\right)
-1\right]  rs^{r-1}t \label{f}%
\end{equation}
or, using (\ref{rays}),%
\begin{equation}
f=x^{r}-\left(  x-t\right)  ^{r}+n\left[  \ln\left(  \frac{n}{t}\right)
-1\right]  . \label{f(x,n,t)}%
\end{equation}

To solve the transport equation (\ref{transport}), we need to compute
$\frac{\partial^{2}f}{\partial n^{2}},\frac{\partial g}{\partial n}$ and
$\frac{\partial g}{\partial x}$ as functions of $t$ and $s.$ Use of the chain
rule gives%
\[%
\begin{bmatrix}
\frac{\partial x}{\partial t} & \frac{\partial x}{\partial s}\\
\frac{\partial n}{\partial t} & \frac{\partial n}{\partial s}%
\end{bmatrix}%
\begin{bmatrix}
\frac{\partial t}{\partial x} & \frac{\partial t}{\partial n}\\
\frac{\partial s}{\partial x} & \frac{\partial s}{\partial n}%
\end{bmatrix}
=%
\begin{bmatrix}
1 & 0\\
0 & 1
\end{bmatrix}
\]
and hence,%
\begin{equation}%
\begin{bmatrix}
\frac{\partial t}{\partial x} & \frac{\partial t}{\partial n}\\
\frac{\partial s}{\partial x} & \frac{\partial s}{\partial n}%
\end{bmatrix}
=\frac{1}{J(t,s)}%
\begin{bmatrix}
\frac{\partial n}{\partial s} & -\frac{\partial x}{\partial s}\\
-\frac{\partial n}{\partial t} & \frac{\partial x}{\partial t}%
\end{bmatrix}
,\label{inversion}%
\end{equation}
where the Jacobian $J(t,s)$ was defined in (\ref{J}). Using (\ref{rays}) and
(\ref{J}) in (\ref{inversion}) we have%
\begin{equation}%
\begin{bmatrix}
\frac{\partial t}{\partial x} & \frac{\partial t}{\partial n}\\
\frac{\partial s}{\partial x} & \frac{\partial s}{\partial n}%
\end{bmatrix}
=\frac{1}{(r-1)t-s}%
\begin{bmatrix}
\left(  r-1\right)  t & -\frac{s^{2-r}}{r}\\
-s & \frac{s^{2-r}}{r}%
\end{bmatrix}
.\label{inversion1}%
\end{equation}
From (\ref{f}) and (\ref{inversion1}) we get%
\begin{align}
\frac{\partial^{2}f}{\partial n^{2}} &  =\frac{\partial^{2}f}{\partial
n\partial t}\frac{\partial t}{\partial n}+\frac{\partial^{2}f}{\partial
n\partial s}\frac{\partial s}{\partial n}=\frac{r-1}{rs^{r-1}\left[
(r-1)t-s\right]  },\nonumber\\
\frac{\partial g}{\partial n} &  =\frac{\partial g}{\partial t}\frac{\partial
t}{\partial n}+\frac{\partial g}{\partial s}\frac{\partial s}{\partial
n}=\frac{\frac{\partial g}{\partial s}-\frac{\partial g}{\partial t}}%
{rs^{r-2}\left[  (r-1)t-s\right]  },\label{fnngxgn}\\
\frac{\partial g}{\partial x} &  =\frac{\partial g}{\partial t}\frac{\partial
t}{\partial x}+\frac{\partial g}{\partial s}\frac{\partial s}{\partial
x}=\frac{\left(  r-1\right)  t\frac{\partial g}{\partial t}-s\frac{\partial
g}{\partial s}}{(r-1)t-s}.\nonumber
\end{align}

Using (\ref{fnngxgn}) in (\ref{transport}), we obtain the ODE%
\begin{equation}
\frac{\partial g}{\partial t}=-\frac{r-1}{2\left[  (r-1)t-s\right]  },
\label{gt}%
\end{equation}
and (\ref{initial}) gives $g(0,s)=0.$ $\allowbreak$Solving (\ref{gt}) we get%
\begin{equation}
g(t,s)=\frac{1}{2}\ln\left[  \frac{s}{s-(r-1)t}\right]  \label{g}%
\end{equation}
or, using (\ref{rays}),%
\begin{equation}
g=\frac{1}{2}\ln\left(  \frac{x-t}{x-rt}\right)  \label{g(x,n,t)}%
\end{equation}
Note that $g$ is undefined when $x=rt,$ i.e., for $x=X_{c}(n).$

\strut Thus, for $\ x>X_{c}(n)$ we have%
\begin{equation}
H_{n}^{r}(x)\sim\exp\left[  x^{r}-\left(  x-t\right)  ^{r}+n\ln\left(
\frac{n}{t}\right)  -n\right]  \sqrt{\frac{x-t}{x-rt}},\quad n\rightarrow
\infty\label{H}%
\end{equation}
with $t(x,n)$ defined implicitly by%
\begin{equation}
r\left(  x-t\right)  ^{r-1}t-n=0. \label{eqt}%
\end{equation}

\subsection{The function $t(x,n)$}

To solve (\ref{eqt}), we shall use Lagrange's inversion formula.

\begin{theorem}
\label{lagrange}Let $\psi(u)$ and $\phi(u)$ be formal power series in $u,$
with $\phi\left(  0\right)  =1.$ Then there is a unique formal power series
$u=u(z)$ that satisfies
\begin{equation}
u=z\phi(u). \label{equ}%
\end{equation}
Further, we have%
\begin{equation}
\left[  z^{k}\right]  \left\{  \psi\left[  u(z)\right]  \right\}  =\frac{1}%
{k}\left[  u^{k-1}\right]  \left\{  \psi^{\prime}(u)\phi(u)^{k}\right\}  ,
\label{equ1}%
\end{equation}
where by $\left[  z^{k}\right]  \left\{  \psi(z)\right\}  $ we mean the
coefficient of $z^{k}$ in the power series of $\psi(z).$
\end{theorem}

\begin{proof}
See \cite[Theorem 5.1]{MR2172781}.
\end{proof}

We first rearrange (\ref{eqt}) so that it looks like (\ref{equ}) and obtain%
\[
\frac{t}{x}=\frac{n}{rx^{r}}\left(  1-\frac{t}{x}\right)  ^{1-r}.
\]
Using (\ref{equ1}) we then have%
\begin{gather*}
\left[  \left(  \frac{n}{rx^{r}}\right)  ^{k}\right]  \left\{  \frac{t}%
{x}\right\}  =\frac{1}{k}\left[  u^{k-1}\right]  \left\{  \left(  1-\frac
{t}{x}\right)  ^{\left(  1-r\right)  k}\right\} \\
=\frac{1}{k}\left[  \left(  \frac{t}{x}\right)  ^{k-1}\right]
{\displaystyle\sum\limits_{j=0}^{\infty}}
\binom{\left(  1-r\right)  k}{j}\left(  -1\right)  ^{j}\left(  \frac{t}%
{x}\right)  ^{j}\\
=\frac{1}{k}\binom{\left(  1-r\right)  k}{k-1}\left(  -1\right)  ^{k-1}.
\end{gather*}
Thus,%
\[
\frac{t}{x}=%
{\displaystyle\sum\limits_{k=1}^{\infty}}
\frac{1}{k}\binom{\left(  1-r\right)  k}{k-1}\left(  -1\right)  ^{k-1}\left(
\frac{n}{rx^{r}}\right)  ^{k}%
\]
or%
\begin{equation}
t_{\mathrm{out}}(x,n)=x%
{\displaystyle\sum\limits_{k=1}^{\infty}}
\frac{1}{k}\binom{\left(  1-r\right)  k}{k-1}\left(  -1\right)  ^{k-1}\left(
\frac{n}{rx^{r}}\right)  ^{k}=\lambda x\left[  1-\rho\left(  \frac{n}{rx^{r}%
}\right)  \right]  , \label{tout}%
\end{equation}
where $\lambda$ was defined in (\ref{lambda}) and
\begin{equation}
\rho\left(  z\right)  =%
{\displaystyle\sum\limits_{k=0}^{\infty}}
\dbinom{rk}{k}\frac{1}{1-rk}z^{k}. \label{rho}%
\end{equation}

To find the region of the complex plane where $\rho\left(  z\right)  $ is
analytic, we compute its radius of convergence using the ratio test. We have%
\[
\underset{k\rightarrow\infty}{\lim}\frac{\dbinom{rk}{k}\frac{1}{1-rk}}%
{\dbinom{r\left(  k+1\right)  }{k+1}\frac{1}{1-r\left(  k+1\right)  }}%
=\frac{\lambda^{r}}{r-1},
\]
where we have used Stirling's formula \cite{stirling} and (\ref{lambda})%
\[
\Gamma(z)\sim\sqrt{\frac{2\pi}{z}}z^{z}e^{-z},\quad z\rightarrow\infty.
\]
Thus, we conclude that $t_{\mathrm{out}}(x,n)$ is analytic for%
\[
\left\vert \frac{n}{rx^{r}}\right\vert <\frac{\lambda^{r}}{r-1}%
\]
or $\left\vert x\right\vert >X_{c}(n),$ with $X_{c}(n)$ defined in (\ref{Xc}).
\ Note that from (\ref{tout}) we have%
\[
t_{\mathrm{out}}(\omega x,n)=\omega t_{\mathrm{out}}(x,n),\text{\quad
for}\quad\omega^{r}=1.
\]

\begin{proposition}
The function $\rho\left(  z\right)  $ satisfies the following properties:

\begin{enumerate}
\item We can represent $\rho\left(  z\right)  $ as a hypergeometric function%
\begin{equation}
\rho\left(  z\right)  =\ _{r-1}F_{r-2}\left[  \left.
\begin{array}
[c]{c}%
\frac{-1}{r},\frac{1}{r},\frac{2}{r},\ldots,\frac{r-2}{r}\\
\frac{1}{r-1},\frac{2}{r-1},\ldots,\frac{r-2}{r-1}%
\end{array}
\right\vert \frac{\left(  r-1\right)  }{\lambda^{r}}z\right]  ,\quad\left\vert
z\right\vert <\frac{\lambda^{r}}{r-1}. \label{rhoF}%
\end{equation}

\item In particular, for $r=2$, we have%
\[
\rho\left(  z\right)  =\sqrt{1-4z},\quad\left\vert z\right\vert <\frac{1}{4}%
\]
and for $r=3,$
\[
\rho\left(  x,n\right)  =\cos\left[  \frac{2}{3}\arcsin\left(  \frac{3}%
{2}\sqrt{3z}\right)  \right]  ,\quad\left\vert z\right\vert <\frac{4}{27}.
\]

\end{enumerate}
\end{proposition}

\begin{proof}
\begin{enumerate}
\item We have%
\begin{align*}
&  _{r-1}F_{r-2}\left[  \left.
\begin{array}
[c]{c}%
\frac{-1}{r},\frac{1}{r},\frac{2}{r},\ldots,\frac{r-2}{r}\\
\frac{1}{r-1},\frac{2}{r-1},\ldots,\frac{r-2}{r-1}%
\end{array}
\right\vert \frac{\left(  r-1\right)  }{\lambda^{r}}z\right] \\
&  =%
{\displaystyle\sum\limits_{k=0}^{\infty}}
{\displaystyle\prod\limits_{j=0}^{r-3}}
\frac{\left(  \frac{j+1}{r}\right)  _{k}}{\left(  \frac{j+1}{r-1}\right)
_{k}}\frac{\left(  -\frac{1}{r}\right)  _{k}}{\left(  1\right)  _{k}}\left[
\frac{\left(  r-1\right)  }{\lambda^{r}}z\right]  ^{k}%
\end{align*}%
\begin{align*}
&  =%
{\displaystyle\sum\limits_{k=0}^{\infty}}
\frac{%
{\displaystyle\prod\limits_{j=0}^{r-1}}
\left(  \frac{j+1}{r}\right)  _{k}}{%
{\displaystyle\prod\limits_{j=0}^{r-2}}
\left(  \frac{j+1}{r-1}\right)  _{k}}\frac{\left(  -\frac{1}{r}\right)  _{k}%
}{\left(  \frac{r-1}{r}\right)  _{k}\left(  1\right)  _{k}}\left[
\frac{\left(  r-1\right)  }{\lambda^{r}}z\right]  ^{k}\\
&  =%
{\displaystyle\sum\limits_{k=0}^{\infty}}
{\displaystyle\prod\limits_{j=0}^{r-1}}
\frac{\Gamma\left(  \frac{j+1}{r}+k\right)  }{\Gamma\left(  \frac{j+1}%
{r}\right)  }%
{\displaystyle\prod\limits_{j=0}^{r-2}}
\frac{\Gamma\left(  \frac{j+1}{r-1}\right)  }{\Gamma\left(  \frac{j+1}%
{r-1}+k\right)  }\frac{1}{\left(  1-rk\right)  k!}\left[  \frac{\left(
r-1\right)  }{\lambda^{r}}z\right]  ^{k}%
\end{align*}%
\begin{align*}
&  =%
{\displaystyle\sum\limits_{k=0}^{\infty}}
\frac{\Gamma\left(  rk+1\right)  }{\Gamma\left[  \left(  r-1\right)
k+1\right]  }\left[  \frac{\left(  r-1\right)  ^{r-1}}{r^{r}}\right]
^{k}\frac{1}{\left(  1-rk\right)  k!}\left[  \frac{\left(  r-1\right)
}{\lambda^{r}}z\right]  ^{k}\\
&  =%
{\displaystyle\sum\limits_{k=0}^{\infty}}
\dbinom{rk}{k}\left(  \frac{\lambda^{r}}{r-1}\right)  ^{k}\frac{1}{\left(
1-rk\right)  }\left[  \frac{\left(  r-1\right)  }{\lambda^{r}}z\right]
^{k}=\rho\left(  z\right)  .
\end{align*}

\item If $r=2,$ we have%
\begin{equation}
\rho\left(  z\right)  =%
{\displaystyle\sum\limits_{k=0}^{\infty}}
\dbinom{2k}{k}\frac{1}{1-2k}z^{k}. \label{rho2}%
\end{equation}
Using the identity \cite[(2.43)]{MR2172781}%
\[%
{\displaystyle\sum\limits_{k=0}^{\infty}}
\dbinom{2k}{k}x^{k}=\left(  1-4x\right)  ^{-\frac{1}{2}},
\]
we have%
\begin{gather*}%
{\displaystyle\sum\limits_{k=0}^{\infty}}
\dbinom{2k}{k}\frac{1}{1-2k}x^{-2k}=\frac{1}{x}\int^{x}%
{\displaystyle\sum\limits_{k=0}^{\infty}}
\dbinom{2k}{k}x^{-2k}\\
=\frac{1}{x}\int^{x}\left(  1-4u^{-2}\right)  ^{-\frac{1}{2}}du=\frac{1}%
{x}\sqrt{x^{2}-4},
\end{gather*}
or%
\[%
{\displaystyle\sum\limits_{k=0}^{\infty}}
\dbinom{2k}{k}\frac{1}{1-2k}x^{k}=\sqrt{1-4x}%
\]
and (\ref{rho2}) follows.

The case $r=3$ was computed using (\ref{rhoF}) and Maple 10.
\end{enumerate}
\end{proof}

We shall now find a representation for $t(x,n)$ in the disc $\left\vert
x\right\vert <X_{c}(n).$ We first observe that when $x=0,$ we have from
(\ref{eqt})%
\begin{equation}
r\left(  -1\right)  ^{r-1}\left[  t(0,n)\right]  ^{r}-n=0, \label{eqt0}%
\end{equation}
which implies that $t(0,n)=\tau_{l}(n),$ where%
\begin{equation}
\tau_{l}(n)=\left[  n\left(  1-\lambda\right)  \right]  ^{1-\lambda}%
\exp\left[  \lambda\left(  2l+1\right)  \pi\mathrm{i}\right]  ,\quad1\leq
l\leq r. \label{t0}%
\end{equation}
Using (\ref{eqt0}) in (\ref{eqt}), we get%
\[
\frac{\left(  t-x\right)  ^{r-1}t}{\left(  \tau_{l}\right)  ^{r}}=1,
\]
or
\begin{equation}
\frac{t}{\tau_{l}}\left(  \frac{t}{\tau_{l}}-\frac{x}{\tau_{l}}\right)
^{r-1}=1. \label{eqt1}%
\end{equation}
To solve (\ref{eqt1}) we use the following Lemma.

\begin{lemma}
Given an algebraic equation of the form%
\begin{equation}
a\left(  a-b\right)  ^{c-1}=1,\quad c\neq0, \label{eq1}%
\end{equation}
we formally have%
\[
a=%
{\displaystyle\sum\limits_{k=0}^{\infty}}
\frac{1}{1-k}\binom{\frac{k-1}{c}}{k}b^{k}.
\]

\end{lemma}

\begin{proof}
Solving for $b$ in (\ref{eq1}) we get%
\[
b=a-a^{\frac{1}{1-c}}.
\]
Letting $\xi=a^{\frac{c}{c-1}}-1,$ we have $a=\left(  \xi+1\right)
^{\frac{c-1}{c}}$ and therefore,%
\[
b=\left(  \xi+1\right)  ^{\frac{c-1}{c}}-\left(  \xi+1\right)  ^{-\frac{1}{c}%
}=\xi\left(  \xi+1\right)  ^{-\frac{1}{c}}%
\]
or%
\begin{equation}
\xi=b\left(  \xi+1\right)  ^{\frac{1}{c}}.\label{eq2}%
\end{equation}
Applying Theorem \ref{lagrange} to (\ref{eq2}) with $\phi\left(  \xi\right)
=\left(  \xi+1\right)  ^{\frac{1}{c}}$ and $\psi\left(  \xi\right)  =\left(
\xi+1\right)  ^{\frac{c-1}{c}},$ we obtain%
\begin{gather*}
\left[  b^{k}\right]  \left\{  a\left(  b\right)  \right\}  =\left[
b^{k}\right]  \left\{  \psi\left[  \xi\left(  b\right)  \right]  \right\}
=\frac{1}{k}\left[  \xi^{k-1}\right]  \left\{  \frac{c-1}{c}\left(
\xi+1\right)  ^{-\frac{1}{c}}\left(  \xi+1\right)  ^{\frac{k}{c}}\right\}  \\
=\frac{c-1}{c}\frac{1}{k}\left[  \xi^{k-1}\right]  \left\{  \left(
\xi+1\right)  ^{\frac{k-1}{c}}\right\}  =\frac{c-1}{c}\frac{1}{k}\binom
{\frac{k-1}{c}}{k-1}=\frac{1}{1-k}\binom{\frac{k-1}{c}}{k}%
\end{gather*}
and the result follows. $\allowbreak$
\end{proof}

Thus, applying the Lemma to (\ref{eqt1}) we find that%
\[
\frac{t}{\tau_{l}}=%
{\displaystyle\sum\limits_{j=0}^{\infty}}
\frac{1}{1-j}\binom{\frac{j-1}{r}}{j}\left(  \frac{x}{\tau_{l}}\right)  ^{j},
\]
which we can write as%

\begin{equation}
t_{\mathrm{in}}(x,n;l)=\tau_{l}(n)+\lambda x\mu\left[  \frac{x}{\tau_{l}%
(n)}\right]  ,\quad1\leq l\leq r, \label{tl}%
\end{equation}
with%
\begin{equation}
\mu\left(  z\right)  =%
{\displaystyle\sum\limits_{k=0}^{\infty}}
\dbinom{\frac{k}{r}}{k}\frac{1}{k+1}z^{k}. \label{mu}%
\end{equation}

\strut Applying the ratio test to (\ref{mu}), we get%
\[
\underset{k\rightarrow\infty}{\lim}\frac{\dbinom{\frac{k}{r}}{k}\frac{1}{k+1}%
}{\dbinom{\frac{k+1}{r}}{k+1}\frac{1}{k+2}}=\lambda^{-\lambda}\left(
1-\lambda\right)  ^{\lambda-1}.
\]
From (\ref{t0}) we have%
\[
\left\vert \tau_{l}(n)\right\vert =\left[  n\left(  1-\lambda\right)  \right]
^{1-\lambda}%
\]
and therefore $t_{\mathrm{in}}(x,n;l)$ is analytic in the sector%
\[
\left\vert x\right\vert <\left[  n\left(  1-\lambda\right)  \right]
^{1-\lambda}\lambda^{-\lambda}\left(  1-\lambda\right)  ^{\lambda
-1}=n^{1-\lambda}\lambda^{-\lambda}=X_{c}(n).
\]
The function $\mu\left(  z\right)  $ doesn't have a simple expression in terms
of elementary functions, except for $r=2,$ when we have%
\[
\mu\left(  z\right)  =1+\frac{x}{2+\sqrt{4+x^{2}}}.
\]

\section{Summary and numerical results}

We summarize the results of the previous section in the following theorem.
Although our analysis was done on the positive real axis, we can extend our
results to the whole complex plane, with the exception of the caustic circle
$\left\vert x\right\vert =X_{c}(n)=\lambda^{-\lambda}n^{1-\lambda}.$

\begin{theorem}
Let%
\[
H_{n}^{r}(x)=\left(  -1\right)  ^{n}\exp\left(  x^{r}\right)  \frac{d^{n}%
}{dx^{n}}\exp\left(  -x^{r}\right)  .
\]
Then, as $n\rightarrow\infty,$

\begin{enumerate}
\item For $\ \left\vert x\right\vert >X_{c}(n),$ we have%
\begin{gather*}
H_{n}^{r}(x)\sim H_{\text{out}}(x,n)=\exp\left(  x^{r}\left\{  1-\left[
1-\lambda+\lambda\rho\left(  \frac{n}{rx^{r}}\right)  \right]  ^{r}\right\}
-n\right) \\
\times\left\{  \frac{n}{\lambda x\left[  1-\rho\left(  \frac{n}{rx^{r}%
}\right)  \right]  }\right\}  ^{n}\sqrt{\frac{1-\lambda\left[  1-\rho\left(
\frac{n}{rx^{r}}\right)  \right]  }{1-\left(  r-1\right)  \left[
1-\rho\left(  \frac{n}{rx^{r}}\right)  \right]  }},
\end{gather*}
with $\lambda=\frac{r-1}{r}$ and%
\[
\rho\left(  z\right)  =%
{\displaystyle\sum\limits_{k=0}^{\infty}}
\dbinom{rk}{k}\frac{1}{1-rk}z^{k}.
\]

\item For $\ \left\vert x\right\vert <X_{c}(n),$ the function $t(x,n)$ is
multivalued and therefore we need to add all the different contributions,%
\begin{gather*}
H_{n}^{r}(x)\sim H_{\text{in}}(x,n)=\sum\limits_{l=1}^{r}\exp\left(
x^{r}-\left\{  x-\tau_{l}(n)-\lambda x\mu\left[  \frac{x}{\tau_{l}(n)}\right]
\right\}  ^{r}-n\right) \\
\times\left\{  \frac{n}{\tau_{l}(n)+\lambda x\mu\left[  \frac{x}{\tau_{l}%
(n)}\right]  }\right\}  ^{n}\sqrt{\frac{\tau_{l}(n)-x\left\{  1-\lambda
\mu\left[  \frac{x}{\tau_{l}(n)}\right]  \right\}  }{r\tau_{l}(n)-x\left\{
1-\left(  r-1\right)  \mu\left[  \frac{x}{\tau_{l}(n)}\right]  \right\}  }},
\end{gather*}
with%
\[
\tau_{l}(n)=\left[  n\left(  1-\lambda\right)  \right]  ^{1-\lambda}%
\exp\left[  \lambda\left(  2l+1\right)  \pi\mathrm{i}\right]
\]
and%
\[
\mu\left(  z\right)  =%
{\displaystyle\sum\limits_{k=0}^{\infty}}
\dbinom{\frac{k}{r}}{k}\frac{1}{k+1}z^{k}.
\]

\end{enumerate}
\end{theorem}

\begin{figure}[ptb]
\begin{center}
\rotatebox{270} {\resizebox{!}{5in}{\includegraphics{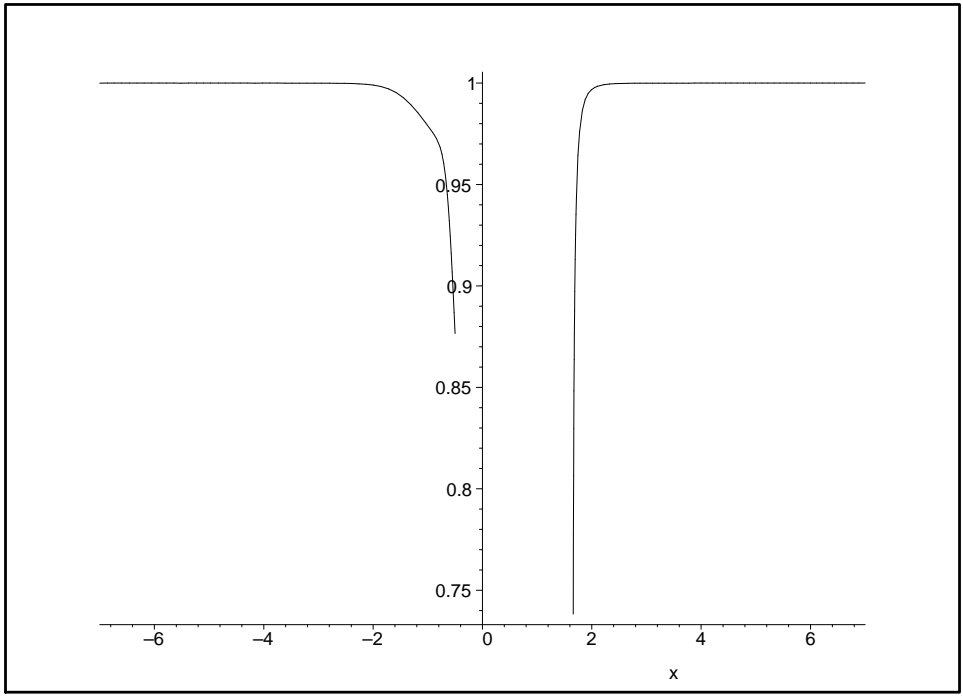}}}
\end{center}
\caption{A sketch of $\frac{H_{5}^{5}(x)}{H_{\text{out}}(x,5)}$ in the outer
region $X_{c}(5) < \left\vert x\right\vert $.}%
\label{out1}%
\end{figure}

In Figure \ref{out1} we sketch the ratio $\frac{H_{n}^{r}(x)}{H_{\text{out}%
}(x,n)},$ for $r=5,$ $n=5$ and $X_{c}(5)\simeq1.649<\left\vert x\right\vert
<7.$ We clearly see how the approximation breaks down inside the caustic
region $\left\vert x\right\vert <X_{c}(n).$

\begin{figure}[ptb]
\begin{center}
$%
\begin{array}
[c]{cccc}%
\rotatebox{270}{\resizebox{2.2in}{!}{\includegraphics{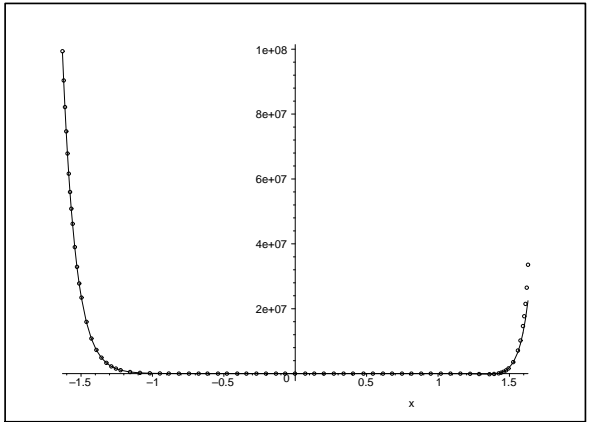}}} &
\rotatebox{270}{\resizebox{2.2in}{!}{\includegraphics{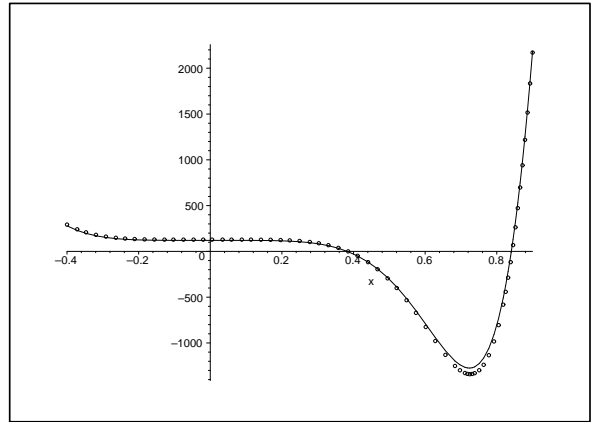}}} &  & \\
\mbox{(a)} & \mbox{(b)} &  & \\
\rotatebox{270}{\resizebox{2.2in}{!}{\includegraphics{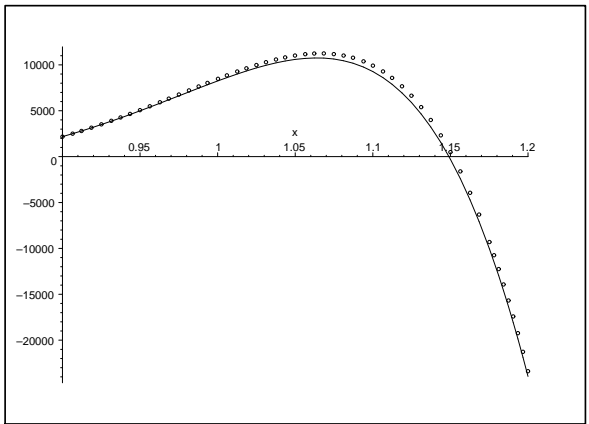}}} &
\rotatebox{270}{\resizebox{2.2in}{!}{\includegraphics{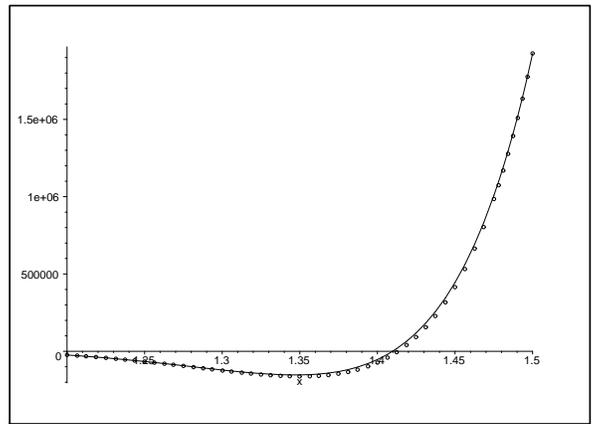}}} &  & \\
\mbox{(c)} & \mbox{(d)} &  &
\end{array}
$
\end{center}
\caption{A comparison of $H_{5}^{5}(x)$ (solid curve) and $H_{\text{in}}(x,5)$
(ooo) inside the caustic region $\left\vert x\right\vert < X_{c}(5)$.}%
\label{in}%
\end{figure}

In Figure \ref{in}(a) we compare the values of $H_{n}^{r}(x)$ and
$H_{\text{in}}(x,n)$ for $r=5,$ $n=5$ and $\left\vert x\right\vert <X_{c}(5).$
To show in detail the graphs for values of $x$ close to the zeros of
$H_{5}^{5}(x)$, we sketch the two functions in smaller intervals in Figures
\ref{in}(b) , \ref{in}(c) and \ref{in}(d).

\begin{figure}[ptb]
\begin{center}
$%
\begin{array}
[c]{cc}%
\rotatebox{270}{\resizebox{2.2in}{!}{\includegraphics{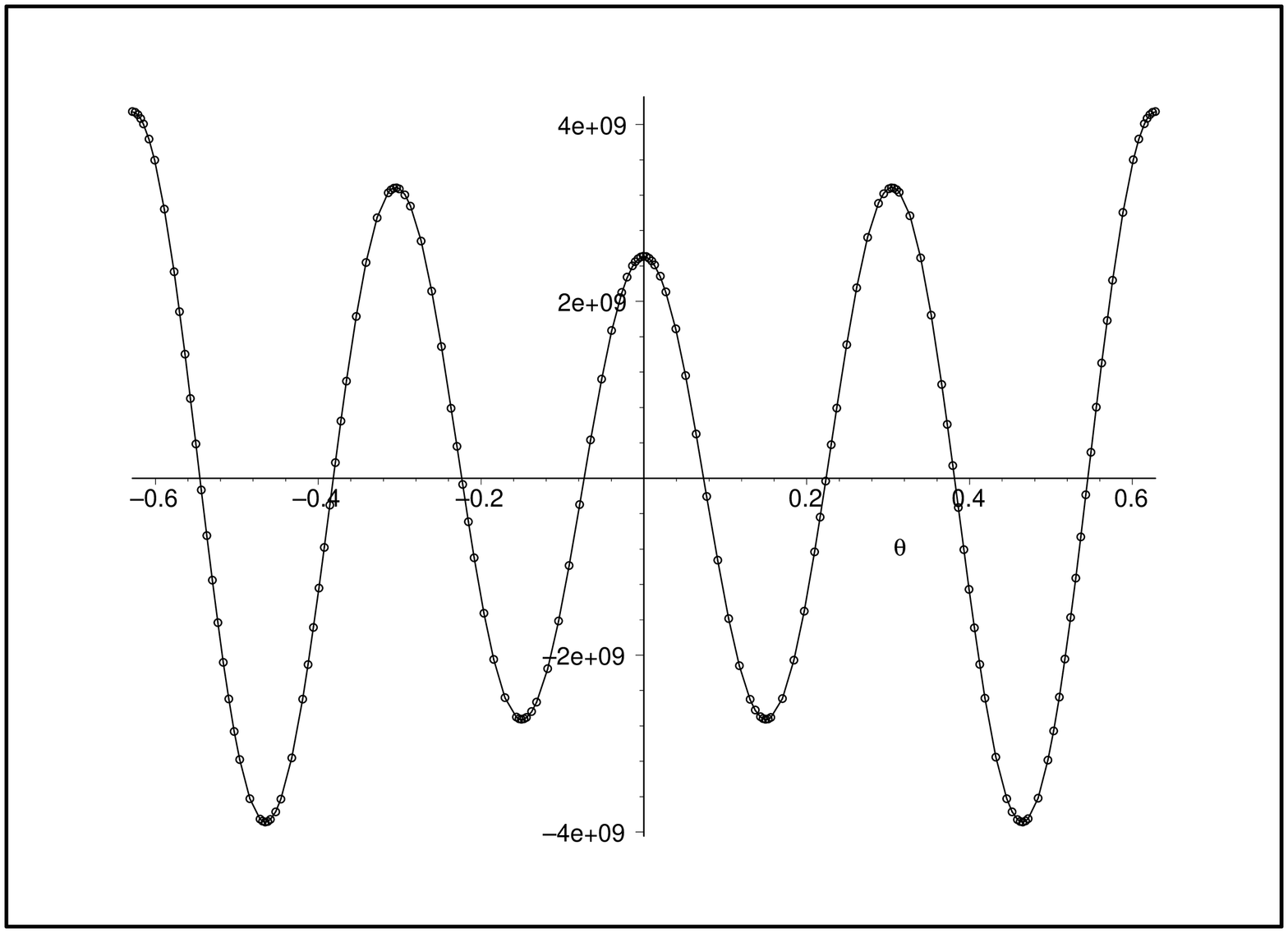}}} &
\rotatebox{270}{\resizebox{2.2in}{!}{\includegraphics{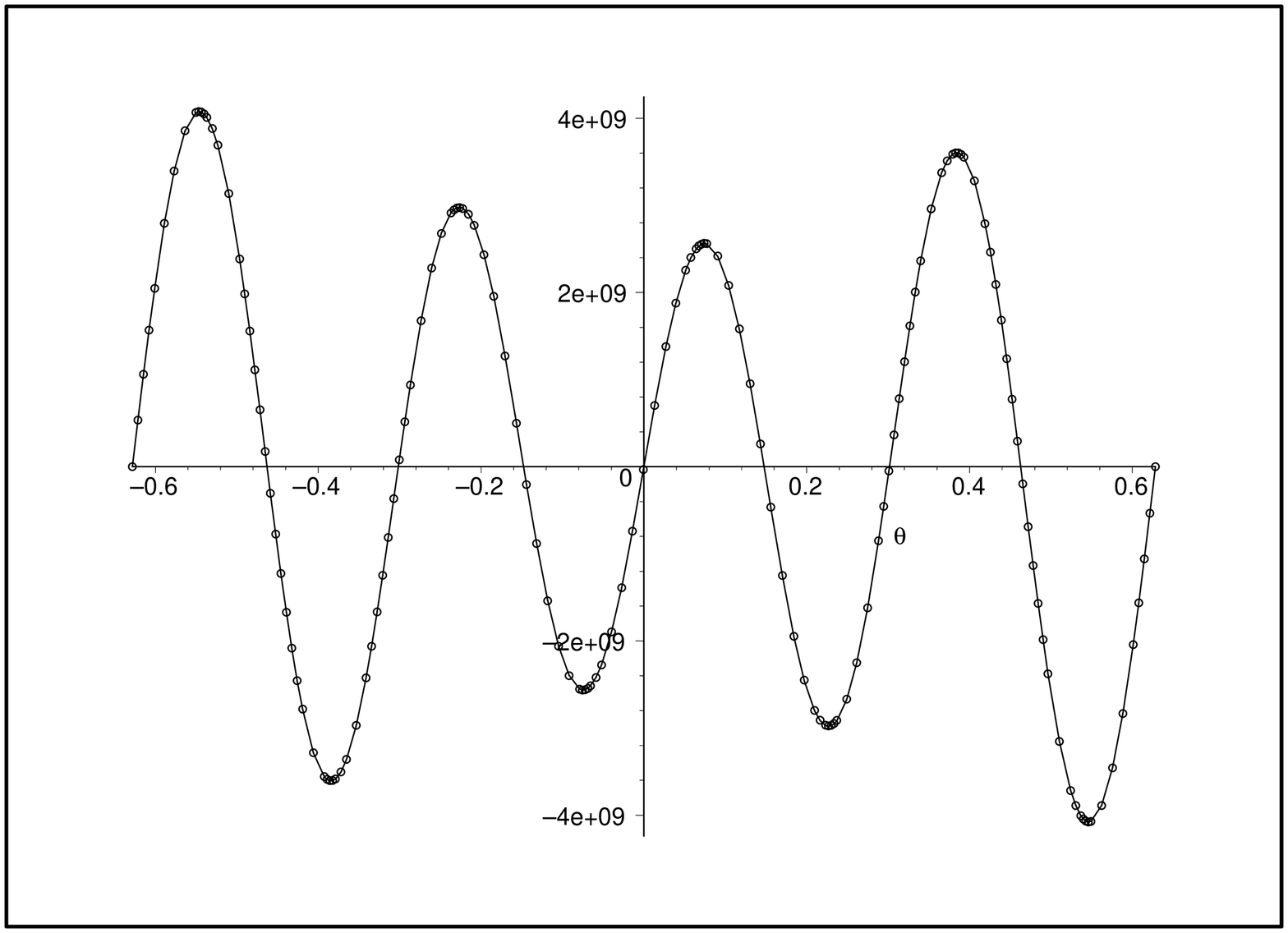}}}\\
\mbox{(a)} & \mbox{(b)}
\end{array}
$
\end{center}
\caption{A comparison of $H_{5}^{5}(x)$ (solid curve) and $H_{\text{out}%
}(x,5)$ (ooo) with $x=2e^{\mathrm{i}\theta}.$}%
\label{circle1}%
\end{figure}

\begin{figure}[ptb]
\begin{center}
$%
\begin{array}
[c]{cc}%
\rotatebox{270}{\resizebox{2.2in}{!}{\includegraphics{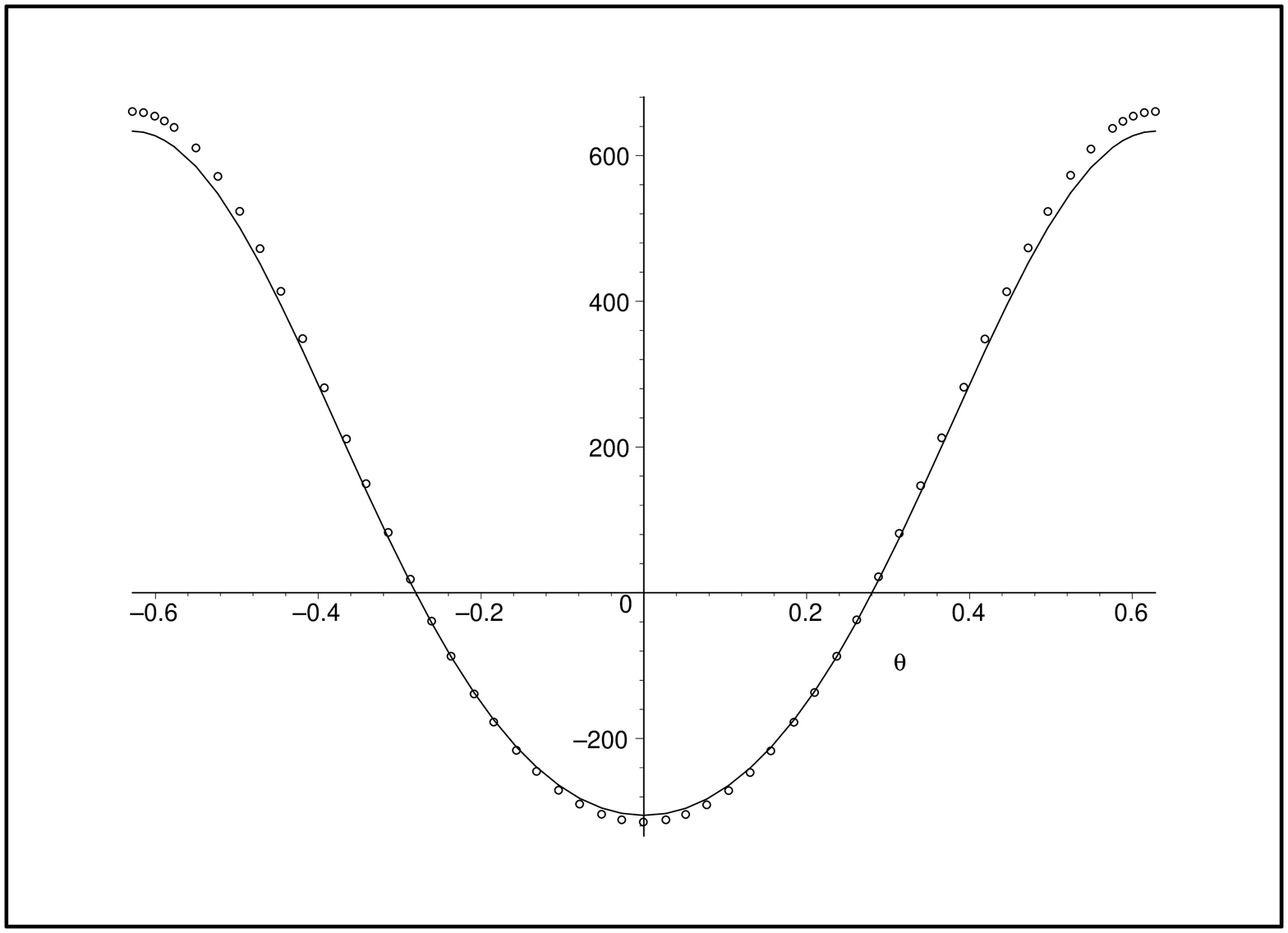}}} &
\rotatebox{270}{\resizebox{2.2in}{!}{\includegraphics{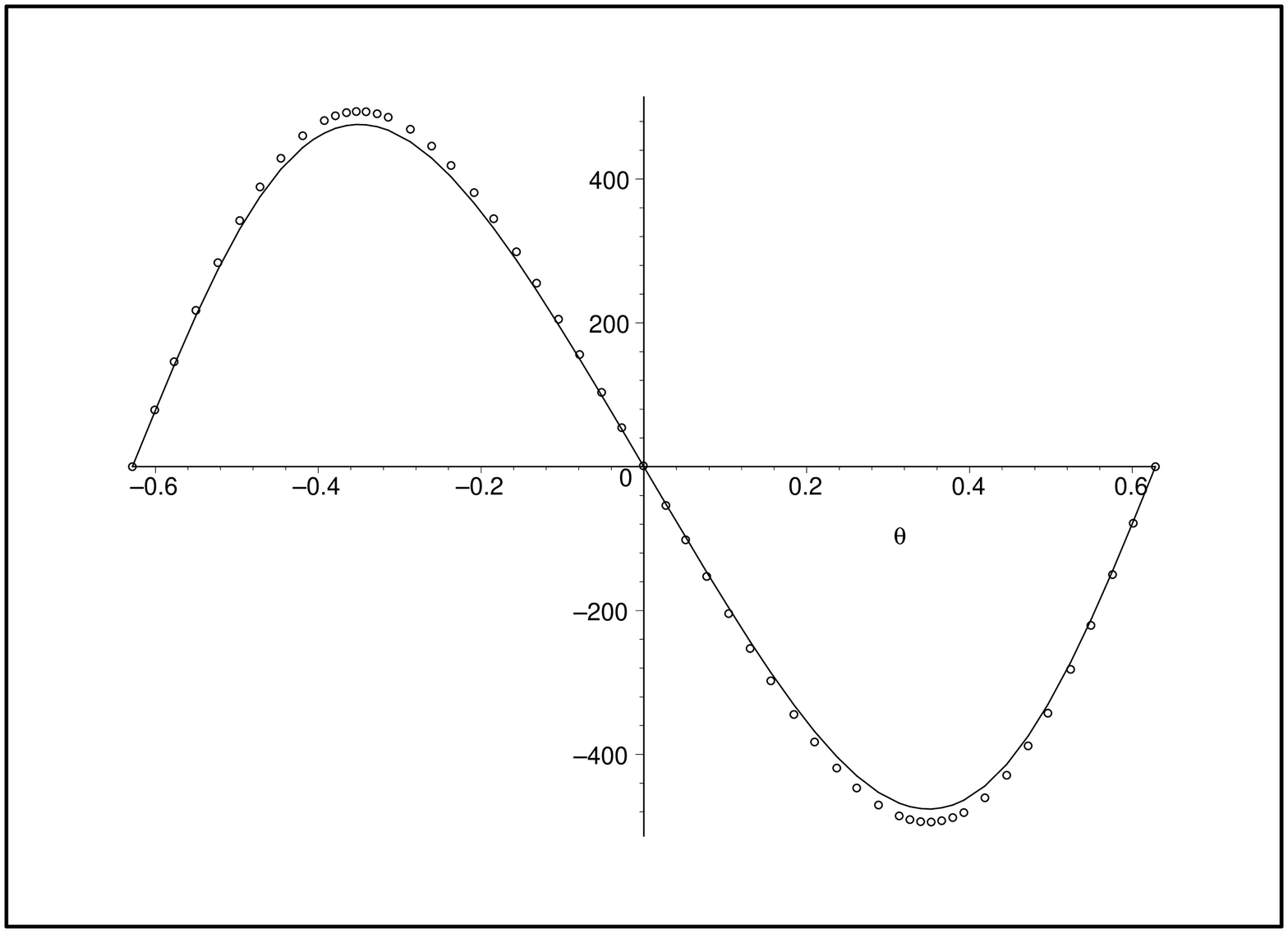}}}\\
\mbox{(a)} & \mbox{(b)}
\end{array}
$
\end{center}
\caption{A comparison of $H_{5}^{5}(x)$ (solid curve) and $H_{\text{in}}(x,5)$
(ooo) with $x=0.5e^{\mathrm{i}\theta}.$}%
\label{circle2}%
\end{figure}

As we mention at the beginning of this section, our approximations are valid
in the complex plane minus the circle $\left\vert x\right\vert =X_{c}(n).$ To
illustrate this, in Figure \ref{circle1} we graph the real (a) and imaginary
(b) parts of $H_{n}^{r}(x)$ and $H_{\text{out}}(x,n)$ for $r=5,$ $n=5$ and
$x=2e^{\mathrm{i}\theta}.$ From Corollary \ref{symm} we know that it is
sufficient to consider the sector $\left\vert \theta\right\vert <\frac{\pi}%
{5}.$ Finally, in Figure \ref{circle2} we do the same for the functions
$H_{n}^{r}(x)$ and $H_{\text{out}}(x,n),$ with $r=5,$ $n=5$ and $x=\frac{1}%
{2}e^{\mathrm{i}\theta}.$

\begin{acknowledgement}
This work was partially supported by a Provost Research Award from SUNY New Paltz.
\end{acknowledgement}

\end{document}